\documentstyle[12pt]{article}

\font\Bbb=msbm10 at 12pt
\font\scriptBbb=msbm10
\font\scriptscriptBbb=msbm7

\newfam\Bbbfam
\textfont\Bbbfam=\Bbb
\scriptfont\Bbbfam=\scriptBbb
\scriptscriptfont\Bbbfam=\scriptscriptBbb

\def\CC{{\fam=\Bbbfam C}}

\def\NN{{\fam=\Bbbfam N}}

\def\RR{{\fam=\Bbbfam R}}

\def\1{1\kern-2.5pt {\rm l}}    

\newcommand{\klc}{\raisebox{-1.5ex}{$\hspace{0.5em}\stackrel{\stackrel{<}
                                          {\textstyle\sim}}{c}\hspace{0.5em}$}}

\newcommand{\schl}[1]{\raisebox{-0.5ex}
                       {$\;\,\stackrel{\textstyle\sim}{#1}\;\,$}}

\newcommand{\gschl}{\raisebox{-1ex}{$\;\stackrel{>}{\sim}\;$}}

\newcommand{\dline}{\||}

\newcommand{\sgn}{$sgn$}

\newcounter{romzahl}

\begin{document}

\begin{Large}
\begin{bf}
\centerline{Vector-valued Lagrange interpolation and}
\centerline{mean convergence of Hermite series}
\end{bf}
\end{Large}         

{\large\centerline{Hermann K\"onig (Kiel)\footnote{Supported by
                   the Brazilian-German CNP$_q$ - GMD agreement}}}

\hspace*{1mm} \\

\underline{Abstract:} 
\begin{enumerate}
\item[]{Let X be a Banach space and $1\le p<\infty$.
We prove interpolation inequalities of Marcinkiewicz-Zygmund type for
X-valued polynomials g of degree $\le n$ on $\RR$,\\
\[
c_p (\sum\limits_{i=1}^{n+1} \mu_i \| g(t_i)e^{-t_i^2 /2} \|^p)^{1/p} \le
(\int\limits_{\RR}^{} \|g(t)e^{-t^2 /2} \|^p dt)^{1/p} \le
d_p (\sum\limits_{i=1}^{n+1} \mu_i \|g(t_i)e^{-t_i^2 /2} \|^p)^{1/p}\;\;,\]
where $(t_i)_1^{n+1}$ are the zeros of the Hermite polynomial $H_{n+1}$
and $(\mu_i)_1^{n+1}$ are suitable weights. The validity of the right
inequality requires $1<p<4$ and X being a UMD-space. This implies a mean
convergence theorem for the Lagrange interpolation polynomials of continuous
functions on $\RR$ taken at the zeros of the Hermite polynomials. In the
scalar case, this improves a result of Nevai $[$N$]$. Moreover, we give
vector-valued extensions of the mean convergence results of Askey-Wainger 
$[$AW$]$ in the case of Hermite expansions.}\\
\end{enumerate}

\section{Introduction and results}

Marcinkiewicz and Zygmund $[$Z, chap.\setcounter{romzahl}{10}\Roman{romzahl}$]$
proved
interpolation inequalities  for trigonometric polynomials g of degree n of
the form
\[
1/3 (\sum\limits_{j=1}^{2n+1} |g(x_j)|^p/(2n+1))^{1/p} \le
(\int\limits_0^{2\pi} |g(x)|^p dx)^{1/p} \le
c_p(\sum\limits_{j=1}^{2n+1} |g(x_j)|^p /(2n+1))^{1/p}
\]
where $x_j=\frac{\pi j}{2n+1}  ,  1<p<\infty$ and $c_p>0$ depends on $p$ only.
The left inequality is true for $p=1,\infty$ as well. A similar result holds
for ordinary polynomials with weight function $(1-x^2)^{-1/2}$, the $x_j$'s
being replaced by the zeros of the Tchebychev polynomials, or, more
generally, for Jacobi polynomial weights $(1-t)^{\alpha}(1+t)^{\beta}$ and
corresponding zeros under suitable restrictions on $p$, provided that
$1/(2n+1)$ is replaced by the weight sequence of the corresponding Gaussian
quadrature formula, which is natural if one considers the case $p=2$.
See Askey $[$A$]$, K\"onig-Nielsen $[$KN$]$.\\

We prove an analogue of this type of inequality in the case of the Hermite
polynomials $H_n$, orthogonal with respect to $e^{-t^2} dt$ on $\RR$.
Let
\begin{enumerate}
\item[]{$h_n := \pi^{-1/4} (2^{n}n!)^{-1/2} H_n ,\;\; {\cal H}_n(t)
 := h_n(t) e^{-t^{2}/2}$}
\end{enumerate}
denote the $L_2$-normalized Hermite polynomials and Hermite functions,
respectively. Thus
\begin{enumerate}
\item[]{$\int\limits_{\RR} h_n(t)h_n(t)e^{-t^2} dt =
\int\limits_{\RR} {\cal H}_n(t) {\cal H}_m(t) dt$,}
\end{enumerate}
see Szeg\"o $[$Sz,5.5$]$.
Let $t_1 >\ldots >t_{n+1}$ denote the zeros of $H_{n+1}$ and
${\lambda}_1,\ldots, {\lambda}_{n+1}$ the corresponding Gaussian quadrature
weights. Thus, for any polynomial q of degree $\le 2n+1$
\begin{enumerate}
\item[]{$\int\limits_{\RR} q(t)e^{-t^2} dt =
 \sum\limits_{j=1}^{n+1} {\lambda}_j q(t_j)$.\hfill (1)}
\end{enumerate}
The Lagrange functions $l_j$, $l_j(t) := h_{n+1}(t)/(h_{n+1}'(t_j)(t-t_j))$,
satisfy $l_j(t_i) = {\delta}_{ij}$. The weights ${\lambda}_j$ can be
calculated from
\begin{enumerate}
\item[]{${\lambda}_j = 2/h_{n+1}'(t_j)^2 = 1/(n h_n(t_j)^2) =
 \int\limits_{\RR} |l_j(t)|^2e^{-t^2} dt$\hfill (2)}
\end{enumerate}
$j = 1,\ldots ,n+1$,  see Szeg\"o $[$Sz, chap. 3.4, 5.5, 15.3$]$. Let
${\mu}_j := {\lambda}_j e^{t_j^2}$. If necessary, we indicate the dependence
on $h_{n+1}$ by superscript, $t_j = t_j^{n+1}, {\lambda}_j = 
{\lambda}_j^{n+1}, {\mu}_j = {\mu}_j^{n+1}$. Recall that
$t_1 \le \sqrt{2n+3}$.\\

Let X be a Banach space and $L_p({\RR};X)$ denote the space of (classes
of) p-th power Boch-\\
ner-integrable functions $f: {\RR} \longrightarrow  X$
with norm $\|f\|_p = (\int\limits_{\RR} \|f(t)\|^p dt)^{1/p}$. Choose
$1<p<\infty$. A Banach space X is an \underline{UMD-space} provided that
the Hilbert-transform on $\RR$
\begin{quote}$Hf(t) := p.v. \int\limits_{\RR} \frac{f(s)}{t-s} ds,
f\in L_p({\RR} ;X)$\end{quote}
defines a bounded linear operator 
$H: L_p({\RR} ;X) \longrightarrow  L_p({\RR} ;X)$. It is well-known that
this holds for \underline{some} $1<p<\infty$ if and only if it holds for
\underline{all} $1<p<\infty$, and the property is thus independent of
$1<p<\infty$, see J.Schwartz $[$S$]$. All $L_q(\mu )$-spaces where
$1<q<\infty$ or all reflexive Orlicz spaces are UMD-spaces, see
Fernandez-Garcia $[$FG$]$.\\

Given $n\in {\NN}$ and a Banach space $X$, the $X$-valued polynomials of
degree $\le n$ will be denoted by 
$\Pi _n(X) = \{\sum\limits_{j=0}^{n} x_jt^j | x_j\in X\}$. We let
$\Pi _n = \Pi _n({\cal K} )$ for ${\cal K} \in
\{{\RR} , {\CC} \}$. Recall that
$(t_j)$ were the zeros of $H_{n+1}$. The Marcinkiewicz--Zygmund interpolation
inequalities for the Hermite polynomials then state:\\

\underline{Theorem 1.} Let X be a Banach space and $1\le p\le\infty$.
\begin{enumerate}
\item[(a)]{The following are equivalent:}
 \begin{enumerate}
 \item[(1)]{There is $c_p>0$ such that for all $n\in {\NN}$ and for all
  $q\in \Pi _n(X)\\[0.25ex]
  (\int\limits_{{\RR}} \|q(t)e^{-t^2/2}\|^p dt)^{1/p} \le
  c_p (\sum\limits_{j=1}^{n+1} \mu _j \|q(t_j)e^{-t_j^2/2}\|^p)^{1/p}$
  \hfill (3)}\\[-3ex]
 \item[(2)]{$X$ is a UMD-space and $1<p<4$.}
 \end{enumerate}
\item[(b)]{}
 \begin{enumerate}
 \item[(1)]{Let $0<\delta <1$. Then there is $c_{\delta } >0$ such that for all
  $n \in {\NN}$ , $N := 2n+3$ and $q\in \Pi _m(X)$ where $m\le 2n$,\\[0.25ex]
  $(\sum\limits_{|t_j|\le \delta \sqrt{N} }\mu _j
  \|q(t_j)e^{-t_j^2/2}\|^p)^{1/p}\le c_{\delta}(\int\limits_{\RR}
                      \|q(t)e^{-t^2/2}\|^p dt)^{1/p}$\hfill (4)}\\[-2.5ex]
 \item[(2)]{If X is a UMD-space and $4/3<p\le \infty$, there is $c_p>0$ such
  that for all n and $q\in \Pi_n(X)\\[0.5ex]
  (\sum\limits_{j=1}^{n+1} \mu _j \|q(t_j)e^{-t_j^2/2}\|^p)^{1/p} \le
  c_p (\int\limits_{\RR} \|q(t)e^{-t^2/2}\|^p dt)^{1/p}$.\hfill (5)\\}\\
 \end{enumerate}
\end{enumerate}
\vspace*{-4ex}\underline{Remarks}
\begin{enumerate}
\item[i)]{As indicated, the $\mu _j$'s and $t_j$'s depend on $n$ as well. It
 seems likely that (5) holds for all $1\le p\le\infty$ and all Banach
 spaces as in the Jacobi case; in (4) the terms involving the zeros
 with largest absolute value are omitted on the left.}
\item[ii)]{The papers of Pollard $[$P$]$, Askey-Wainger $[$AW$]$ and Nevai
 $[$N$]$ strongly suggest the choice of the weight function $e^{-p/2 t^2}$
 instead of $e^{-t^2}$ for $p\not= 2$ if positive results are to be expected.}
\end{enumerate}
Inequalities of type (3) imply mean convergence results for interpolating
polynomials. Given a continuous function $f: {\RR} \longrightarrow X,\;\;
I_nf := \sum\limits_{j=1}^{n+1} f(t_j)l_j$
is the interpolating polynomial at the
zeros of the Hermite polynomial $H_{n+1},\;\; I_nf(t_j) = f(t_j)$.\\
Let $L_p({\RR} , e^{-p/2 t^2};X) =
\{ f: {\RR} \longrightarrow X | g\in L_p ({\RR} ;X)$
where $g(t) = f(t)e^{-t^2/2}\}$.\\

\underline{Theorem 2.}
\begin{enumerate}
\item[]{Let X be a UMD-space and $1<p<4$. Let $\alpha >1/p$ and
 $f: {\RR} \longrightarrow X$ be be a continuous function satisfying
  \begin{enumerate}
  \item[]{$\lim\limits_{|t|\to\infty} \|f(t)\|_X (1+|t|)^{\alpha}
                                       e^{-t^2/2} = 0$.\hfill (6)}
  \end{enumerate}
 Then the interpolating polynomials $I_nf$ at the zeros of the Hermite
 polynomials $H_{n+1}$ converge to f in the
 $L_p({\RR} , e^{-p/2t^2};X)$-norm,
  \begin{enumerate}
  \item[]{$(\int\limits_{{\RR}} \|(f(t)-I_nf(t)) e^{-t^2/2}\|^p dt)^{1/p}
   \longrightarrow  0$\hfill (7)}
  \end{enumerate}
 The same statement is false, in general, if $p>4$.\\}
\end{enumerate}

\underline{Remarks}
\begin{enumerate}
\vspace{-0.5ex}
\item[(i)]{In the scalar case $X={\cal K}$, this improves the result of Nevai
 $[$N$]$ where $\alpha =1$ is assumed. Nevai's result, however, holds
 for any $1<p<\infty$.}
\vspace{-0.5ex}
\item[(ii)]{Condition (5) clearly guarantees that
 $f\in L_p({\RR} , e^{-p/2t^2};X)$. Under the weaker assumption
 that $f\in C({\RR} ;X)\cap L_p({\RR} ; e^{-p/2t^2};X)$
 convergence (7) does not hold, in general, as we show below.}
\end{enumerate}

Askey-Wainger $[$AW$]$ prove their mean convergence result for the expansions
of the $L_p$-functions on ${\RR}$ into Hermite functions for $4/3<p<4$.
Their proof generalizes to the vector-valued setting if X is a UMD-space.
The necessity of the UMD condition of the following result is a consequence
of theorem 1. Given $f\in L_p({\RR} ;X)$, we let
$a_j = \int\limits_{{\RR}} f(t){\cal H}_j(t) dt$,
$P_nf = \sum\limits_{j=0}^{n} a_j{\cal H}_j$. We have:\\

\underline{Theorem 3.}
\begin{enumerate}
\item[]{Let X be a Banach space and $1\le p\le\infty$. The following are
 equivalent:}
  \begin{enumerate}
  \item[(1)]{For all $f\in L_p({\RR} ;X)$, $P_nf\longrightarrow f$ in
   $L_p$-norm, i.e.}
   \begin{enumerate}
   \item[]{$\int\limits_{\RR} \|f(t)-P_nf(t)\|^p dt \longrightarrow 0$.}
   \end{enumerate}
  \item[(2)]{X is a UMD-space and $4/3<p<4$.}
  \end{enumerate}
\end{enumerate}
  
Work on this paper started during a visit of the University of Campinas. The
author greatfully acknowledges the hospitality of the collegues there, in
particular D.L.Fernandez and K.Floret.\\
                                        
\section{Hermite asymptotics, zeros and quadrature weights}

We need some estimates for the quadrature weights and the zeros of the
Hermite polynomials. The zeros $t_1>\ldots >t_{n+1}$ of $H_{n+1}$ are in in the
interval $(-\sqrt{N}, \sqrt{N})$ where $N := 2n+3$; $t_{[\frac{n+2}{2}]}$ is
the one closest to 0 and, more precisely, $\sqrt{N}-t_1 =
{\cal O} (1/\sqrt[6]{N})$, Szeg\"o $[$Sz, 6.32$]$. Given sequences
$(a_n)_{n\in {\NN}}, (b_n)_{n\in {\NN}}$ and $c>0$, we write
$a_n\klc b_n$ if $a_n \le c\,b_n$
holds for all $n\in {\NN}$.
We write $a_n$ \schl{c} $b_n$ if $a_n\klc b_n$ and $b_n\klc a_n$.
A similar notation will be used for real functions.
Define $\Phi : [0,1] \longrightarrow {\RR}$ by
\begin{enumerate}
\item[]{$2/3\;\Phi (x)^{3/2} = \int\limits_x^1 (1-s^2)^{1/2}\,ds ,\;\;
x\in [0,1]$.}
\end{enumerate}
One checks easily that $\Phi (x)$ \schl{c} $(1-x^2)$ with
$c = 2^{2/3}$ and that\\
$|\Phi '(x)|^{-1/2} =
(\Phi (x)\,/\,(1-x^2))^{1/4} \; \schl{d} \; 1$ with $d = 2^{1/6}$.\\
Skovgaard's asymptotic formula for $H_{n+1}$ yields that for
$0 \le t \le \sqrt{N}-\sqrt[1/6]{N}$
\begin{enumerate}
\item[]{${\cal H}_{n+1}(t) = c_n/n^{1/12} \;
 |\Phi '(t/\sqrt{N})|^{-1/2}
 \; \{ Ai \: (-N^{2/3} \, \Phi (\frac{t}{\sqrt{N}}))
 \: + {\cal O}(\frac{1}{n^{7/6} \, \Phi (t/\sqrt{N})^{1/4}})\}$}
\end{enumerate}
where $\lim\limits_n c_n = 2^{5/12}$ , cf. Askey-Wainger $[$AW, p.700$]$
(there is a misprint, it should be $H_n$ instead of ${\cal H}_n$).

\newpage
Using the relation between Airy and Bessel functions
\begin{enumerate}
\item[]{$Ai(-z) = \sqrt{z} /3 \: ($J$_{1/3}(\zeta) + $J$_{-1/3}(\zeta)) =
 \frac{1}{\sqrt{\pi}}\frac{1}{\sqrt[4]{z}} \;
 \cos (\zeta\,-\,\pi /4) \: (1\,+\,{\cal O}(\frac{1}{\zeta}))$}
\end{enumerate}
for large \ $\zeta := 2/3\: z^{3/2}$
\begin{enumerate}
\item[]{${\cal H}_{n+1}(t) = \frac{g_n(t)}{n^{1/8}\,(\sqrt{N}-t)^{1/4}} \:
 \cos (\frac{2}{3}N \, \Phi (\frac{t}{\sqrt{N}})^{3/2} \: -\frac{\pi}{4})
 \; + {\cal O}(\frac{1}{n^{9/8}\,(\sqrt{N}-t)^{1/4}})$\hfill (8)}
\end{enumerate}
where $\; g_n(t) \: \schl{c} \: 1 \;$ for some $c>0$
independent of n and $t\le \sqrt{N}-\sqrt[1/6]{N}$.
We note that $\; N \, \Phi (t/\sqrt{N})^{3/2} \sim n^{1/4}(\sqrt{N}-t)^{3/2}$.\\

Define \ $k_{n,j} := $sup$\{ |{\cal H}_{n+1}(t)| : t\in [t_{j+1},t_j]\}$
for $n\in {\NN} , j = 1,\ldots ,n$ .
By (8), there is $c>0$ independent of $n$ and $j$ such that
\begin{enumerate}
 \item[]{$k_{n,j}$ \schl{c} $n^{-1/8} \, (\sqrt{N}-|t_j|)^{-1/4}$ .\hfill (9)\\}
\end{enumerate}

\underline{Lemma 1.}
\vspace{-2ex}
\begin{enumerate}
\item[]{There is a constant $c>0$ such that for all $n\in{\NN}$ and
 $j = 1,\ldots ,n$\\
 $(t_j-t_{j+1}) \schl{c}  {\mu}_j =
 \frac{2}{n{\cal H}_n(t_j)^2}
  = \frac{1}{
  {\cal H}_{n+1}'(t_j)^2}  
  \schl{c}
  {\frac{1}{n^{1/4}}}{\frac{1}{\sqrt{{\sqrt{N}}-|t_j| }}}$\\
 The sequence $({\mu}_j)_{j=1,\ldots ,[(n+2)/2]}$ is decreasing in j.}
\end{enumerate}

\underline{Proof.}
\begin{enumerate}
\item[]{The equalities follow from (2) and ${\mu}_j = {\lambda}_je^{t_j^2}$.
 By symmetry we may assume $t_j\ge 0$,\\
 $j\le[\frac{n+2}{2}].\;\; {\cal H}_{n+1}$
 satisfies the differential equation}
 \begin{enumerate}
 \item[]{${\cal H}_{n+1}''(t) + (N-t^2){\cal H}_{n+1} = 0,\;\;
  N = 2n+3$,\hfill (10)}
 \end{enumerate}
\end{enumerate}
\begin{enumerate}
\item[]{Szeg\"o $[$Sz, 5.5$]$. Let $f(t) := {\cal H}_{n+1}'(t)^2 +
 (N-t^2){\cal H}_{n+1}(t)^2$.\\
 Then $f'(t) = -2t{\cal H}_{n+1}(t)^2 \le  0$
 for $t\ge 0$, i.e. f is decreasing. Hence ${\mu}_j = 1/{\cal H}_{n+1}'(t_j)^2
 \ge  1/{\cal H}_{n+1}'(t_{j+1})^2 = {\mu}_{j+1}$, using
 ${\cal H}_{n+1}(t_i) = 0$.\\
 Denote
 ${\nu}_{n,j} := \sup\{|{\cal H}_{n+1}'(t)| : t\in [t_{j+1},t_j]\}$. In view
 of (10), ${\cal H}_{n+1}$ is concave or convex in $(t_{j+1},t_j)$,
 depending on whether ${\cal H}_{n+1}$ is positive or negative there. Thus}
 \begin{enumerate}
 \item[]{${\nu}_{n,j} = \max\{ |{\cal H}_{n+1}'(t_{j+1})| ,
  |{\cal H}_{n+1}'(t_j)|\} = |{\cal H}_{n+1}'(t_{j+1})|$.}
 \end{enumerate}
Let $\overline{t}\in (t_{j+1},t_j)$ be such that
${\cal H}_{n+1}'(\overline{t}) = 0$, i.e.
$\kappa_{n,j} = |{\cal H}_{n+1}(\overline{t})|$.
Using that $f$ is decreasing, we find
 \begin{enumerate}
 \item[]{${\nu}_{n,j}^2 = {\cal H}_{n+1}'(t_{j+1})^2 \ge
  (N-{\overline{t}}^2){\kappa}_{n,j}^2 \ge
  {\cal H}_{n+1}'(t_j)^2 = {\nu}_{n,j-1}^2$,}
 \end{enumerate}
and using (9) and $N-{\overline{t}}^2  \schl{}  {\sqrt{n}}({\sqrt{N}}-t_j)$,
 \begin{enumerate}
 \item[]{${\nu}_{n,j}  \schl{c}   {\kappa}_{n,j}\cdot
  {\sqrt[4]{n}}{\sqrt{{\sqrt{N}}-t_j}}  \schl{c}  
  n^{1/8}({\sqrt{N}}-t_j)^{1/4}$.}
 \end{enumerate}
Hence ${\kappa}_{n,j}\cdot{\nu}_{n,j}  {\schl{d}}  1$ with c,d independent of
n and $j = 1,\ldots ,[(n+2)/2]$. Since ${\mu}_{j+1} = 1/{\nu}_{n,j}^2  
(j\ge 2)$, the right estimate
${\mu}_j  \schl{}   n^{-1/4} ({\sqrt{N}}-t_j)^{-1/2}$ follows.

The mean value theorem, applied to ${\cal H}_{n+1}$ in $(t_{j+1},t_j)$,
yields
 \begin{enumerate}
 \item[]{$t_j-t_{j+1} \ge
  \frac{\kappa_ {n,j}}{\nu _{n,j}}  =
  \frac{\kappa_ {n,j}\cdot\nu _{n,j}}{\nu _{n,j}^2}  \schl{d}  
  \frac{1}{\nu _{n,j}^2} = \mu _{j+1}  {\schl{}}  \mu _j$.\hfill (11)}
 \end{enumerate}
For $t\in [-t_j,t_j]$, $N-t^2 \ge  N-t_j^2$. Comparing the differential
equation (10) in the interval $[-t_j,t_j]$ with
\begin{displaymath}
 {\cal K}''(t)+(N-t_j^2){\cal K}(t) = 0, t\in [-t_j,t_j],
\end{displaymath}
Sturm's comparison principle yields the converse to (11),
 \begin{enumerate}
 \item[]{$t_j-t_{j+1} \le  {\frac{\pi}{{\sqrt{N-t_j^2}}}} \le
  \frac{\pi}{\sqrt[4]{n} \sqrt{\sqrt{N}-t_j}} \schl{c}
  {\mu}_j$.\hfill\raisebox{-2ex}{$\Box$}\\}
 \end{enumerate}
\end{enumerate}

More precise information on the constants involved can be found from (8),
analyzing the zeros of the cosine-term. We do not need this. We note however,
that for all $j$ with $|t_j| < \delta\sqrt{N}$ with fixed $0 < \delta  < 1$,
one has ${\mu}_j \schl{c_{\delta}}  n^{-1/2}$, cf. Nevai $[$N$]$. 
In contrast to this, ${\mu}_1  {\schl{c}}   n^{-1/6}$. It is in the range 
in between that lemma 1 is of importance. As a corollary, we get\\

\underline{Lemma 2.}
\begin{enumerate}
\item[]{Let $1\le p\le 2$. Then there is $c>0$ such that for all
 $n\in{\NN}$}
 \begin{enumerate}
 \item[]{$\sum\limits_{j=1}^{n+1}{{\mu}_j^p}\le  c n^{1-p/2}$
 (bounded for $p=2$)} 
 \end{enumerate}
\end{enumerate}

\underline{Proof.}
\begin{enumerate}
\item[]{For \ $p=1$ \ see Nevai $[$N$]$. Let ${\alpha} := p-1$.
 Then by lemma 1
 \[
 {\sum\limits_{j=1}^{n+1}}{\mu}_j^p\schl{c_1}
 {\sum\limits_{j=1}^{n}}(t_j-t_{j+1}){\mu}_j^{\alpha} \schl{c_2}
 {\sum\limits_{j=1}^{n}}(t_j-t_{j+1}){\frac{1}{n^{{\alpha}/4}}}
 {\frac{1}{({\sqrt{N}-|t_j|)^{{\alpha}/2}}}}\;\;\;.\]
 Since $1/{({\sqrt{N}}-t)}$ is monotone for $t\ge 0$ and $(t_j-t_{j+1})$
 is monotone in j (for the positive $t_j$'s), the last Riemann sum can be
 replaced by an integral}
 \begin{enumerate}
 \item[]{$\sum\limits_{j=1}^{n+1}{{\mu}_j^p} \schl{c_3}
  {\frac{1}{n^{{\alpha}/4}}}
  {\int\limits_{-{\sqrt{N}}}^{{\sqrt{N}}}}
  {\frac{dt}{({\sqrt{N}}-|t|)^{{\alpha}/2}}} \schl{c_4}
  n^{(1-{\alpha})/2} = n^{1-p/2}$.\hfill\raisebox{-2ex}{$\Box$}\\}
 \end{enumerate}
\end{enumerate}

\vspace{2ex}
\section{The interpolation inequalities}

To prove theorem 1, we need two well-known facts about continuity in $L_p$, cf.
Benedek-Murphy-Panzone $[$BMP$]$ or Pollard $[$P$]$.\\

\underline{Lemma 3.}
\vspace{-1ex}
\begin{enumerate}
\item[]{Let X be a Banach space, $1\le p\le\infty$, $({\Omega},{\mu})$ be
 a measure space and $k,r: {\Omega}^2\longrightarrow{\cal K}$ be measurable
 such that}
 \begin{enumerate}
 \item[]{$\sup\limits_{u}
  {{\int\limits_{\Omega}^{}|k(u,v)|\;|r(u,v)|^{p'}}\,d{\mu}(v)\le M},
  \sup\limits_{v}
  {\int\limits_{\Omega}^{}
  {|k(u,v)|\;|r(u,v)|^{-p}\,d{\mu}(u)}\le M}$\hfill (12)}
 \end{enumerate}
\item[]{Then $T_kf(u)
 := {\int\limits_{\Omega}^{}{k(u,v)f(v)\,d{\mu}(v)}}$
 defines an operator\\
 $T_k: L_p({\Omega},{\mu};X)\longrightarrow L_p({\Omega},{\mu};X)$ of norm
 $\le M$. Here $1/p+1/p' = 1$.}
\end{enumerate}

This follows from an application of H\"older's inequality. One consequence
is\\

\underline{Lemma 4.}
\vspace{-1ex}
\begin{enumerate}
\item[]{Let X be a Banach space, $1\le p\le\infty, b\in{\RR}$ and
 $k: {\RR}^2\longrightarrow {\RR}$ be defined by\\
 $k(u,v) := |\; |\frac{u}{v}|^b-1|/|u-v|$.
 Then $T_k$ is bounded as a map\\
 $T_k: L_p({\RR};X)\longrightarrow L_p({\RR};X)$
 provided that $-1/p < b < 1-1/p$\\
 (actually if and only if).}
\end{enumerate}

\newpage
For the convenience of the reader, here is a sketch of the proof
(cf. $[$BMP$]$):

Take $r(u,v) := |u/v|^{1/{pp'}}$. To check (12), substitute
$v/u = t$ to find 
\begin{displaymath}
\sup\limits_{u (\not= 0)}\,
\int\limits_{\RR}^{}{|k(u,v)|r(u,v)^{p'}dv} =
\int\limits_{\RR}^{}{|t^{-b}-1|\; |t|^{-1/p}/|t-1|\,dt}.
\end{displaymath}
This is finite: integrability at 0 is assured since $b<1-1/p$, integrability
at $\pm\infty$ since $b>-1/p$. Note that for $t\longrightarrow 1$, the
integrand tends to $|b|$. The second condition in (12) is checked
similarly.\hfill $\Box$\\

Instead of using this simple lemma 4 below one could also apply the
general theory of weighted singular integral operators with weights
in the Muckenhaupt class $A_p$, see.~Garcia-Cuerva,
Rubio de Francia [GR, ch.~\setcounter{romzahl}{4}\Roman{romzahl}].\\

\underline{Proof of theorem 1.}\\
\begin{enumerate}
\item[a)]{$(2)\Rightarrow (1)$}
 \begin{enumerate}
 \item[]{We prove inequality (3) if X is an
  UMD-space and $1<p<4$.
  Let $q\in {\Pi}_n(X)$ and put
  $y_j := {q(t_j)e^{-t_j^2/2}}/{(n^{1/8}{\cal H}_{n+1}'(t_j))}$.
  Then q concides with its interpolating polynomial}
  \begin{enumerate}
  \item[]{$q(t) = I_nq(t) =
   \sum\limits_{j=1}^{n+1}{q(t_j)l_j(t)} =
   \sum\limits_{j=1}^{n+1}{q(t_j){\frac{h_{n+1}(t)}{h_{n+1}'(t_j)(t-t_j)}}}$}
  \end{enumerate}
 \item[]{and we have to estimate}
  \begin{enumerate}
  \item[]{$J :=
   (\int\limits_{{\RR}}^{}{\|q(t)e^{-t^2/2}\|^p dt)^{1/p}} =
   (\int\limits_{{\RR}}^{}{\|n^{1/8}{\cal H}_{n+1}(t)
    \sum\limits_{j=1}^{n+1}{y_j/(t-t_j)\|^p dt)^{1/p}}}$\hfill (13)} 
  \end{enumerate}
 \item[]{from above. Let $I_j := (t_j,t_{j-1}),
  |I_j| = (t_{j-1}-t_j)$ and ${\chi}_j$ be the characteristic function
  of $I_j$, for $j =1,\ldots ,n+1$ (with $t_0 := {\sqrt{N}}, N = 2n+3$).
  The proof uses essentially that $1/(t-t_j)$ is close enough
  to the Hilbert transform of $-{\chi}_j/|I_j|$ at $t$ which is}
  \begin{enumerate}
  \item[]{$H(-\frac{{\chi}_j}{|I_j|}(t) =
   \frac{1}{|I_j|}$ ln$|\frac{t-t_{j-1}}{t-t_j}| =
   \frac{1}{|I_j|}$ ln$|1-\frac{|I_j|}{t-t_j}|$ .}
  \end{enumerate}
 \item[]{We claim that for all $j = 1,\ldots,n+1$}
  \begin{enumerate}
  \item[]{$n^{1/8} |{\cal H}_{n+1}(t)|\;
   |\frac{1}{t-t_j}+H(\frac{{\chi}_j}{|I_j|})(t)|
   \le f_j(t) , t\in{\RR}$ \hfill (14)}
  \end{enumerate}
 \item[]{where}
  \begin{enumerate}
  \item[]{$f_j(t) = c_1$ min$(\frac{1}{|I_j|},
   \frac{|I_j|}{(t-t_j)^2}) |\sqrt{N}-|t|\,|^{-1/4}$ .}
  \end{enumerate}
 \item[]{By Skovgaard's asymptotic formula $[$AW, p.700$]$
  ( (8) for $|t|<\sqrt{N}$ )}
  \begin{enumerate}
  \item[]{$n^{1/8} |{\cal H}_{n+1}(t)| \le 
   c_2 |\sqrt{N}-|t|\; |^{-1/4} , t\in{\RR}$ .}
  \end{enumerate}
 \item[]{Thus, for $|t-t_j|>2|I_j|$ , (14) follows from $|x-$ln$(1+x)| \le x^2$
  for $|x| \le 1/2$ , i.e.}
  \begin{enumerate}
  \item[]{$|\frac{1}{t-t_j}+H(\frac{{\chi}_j}{|I_j|})(t)|
   \le \frac{|I_j|}{(t-t_j)^2}$}
  \end{enumerate}
 \item[]{For $|t-t_j| \le 2|I_j|$ , one uses that ${\cal H}_{n+1}$ has
  zeros at $t_j$ to get (14).
  We remark that for j=1, there is a singularity of the
  $H({\chi}_1/|I_1|)$-term at $t_0=\sqrt{N}$
  (where ${\cal H}_{n+1}$ has no zero to compensate);
  in this case use $|$ln$|x|\; | \le 4|x|^{-1/4}$ for $|x| \le 1$
  and $|I_1| \sim n^{-1/6},\;\;
  |{\cal H}_{n+1}(t)| \stackrel{<}{\sim} n^{-1/12}$ \ \ $[$Sz, 6.32$]$,
  $[$AW$]$ to find}
  \begin{enumerate}
  \item[]{$n^{1/8} |{\cal H}_{n+1}(t)| |H(\frac{{\chi}_1}{|I_1|})(t)|
   \le c_3/(|I_1| |\sqrt{N}-|t|\; |^{1/4})$}
  \end{enumerate}
 \item[]{for $|t-t_1| < 2|I_1|$.
  Now (13) and (14) imply $J \le J_1+J_2$ where}
  \begin{enumerate}
  \item[]{$J_1 = c_2 (\int\limits_{\RR}^{} \|H(\sum\limits_{j=1}^{n+1}
   y_j {\chi}_j/|I_j|)(t)\|^p
   |\sqrt{N}-|t|\; |^{-p/4} dt)^{1/p}$\\
   $J_2 = (\int\limits_{\RR}^{} (\sum\limits_{j=1}^{n+1} \|y_j\|
   f_j(t))^p dt)^{1/p}$ .}
  \end{enumerate}
 \item[]{We estimate the "main term" $J_1$ and the "error-term" $J_2$
  seperately.\\
  
  By lemma 4, the kernel $\frac{1}{|t-s|} | |\frac{s}{t}|^{1/4}-1|$
  defines a bounded operator in $L_p({\RR};X)$ for $1<p<4$.
  Since the Hilbert transform with kernel $1/(t-s)$ is bounded in
  $L_p({\RR};X)$ by assumption, so is the weighted Hilbert transform
  with kernel $|s/t|^{1/4}/(t-s)$.
  Replacing $s$ and $t$ by $\sqrt{N}-s$ and $\sqrt{N}-t$, we find
  that for $g\in L_p({\RR};X)$}
  \begin{enumerate}
  \item[]{$(\int\limits_{\RR}^{}\|\int\limits_{\RR}^{}
   \frac{g(s)}{t-s} |\frac{\sqrt{N}-s}{\sqrt{N}-t}|^{1/4}
   ds \|^p dt)^{1/p} \le c_p (\int\limits_{\RR}^{}
   \|g(s)\|^p ds)^{1/p}$,}
  \end{enumerate}
 \item[]{$c_p$ independent of $g$ and $N$.
  Replacing $t$ by $(-t)$ and $g$ by $g^{-}$,
  $g^{-}(s) = g(-s)$, and then putting $f(s) = |\sqrt{N}-|s|\; |^{1/4} g(s)$,
  we find}
  \begin{enumerate}
  \item[]{\hspace*{-0.5em}$(\int\limits_{\RR}^{} \|\int\limits_{\RR}^{}
   \frac{f(s)}{t-s} ds\|^p\; (\sqrt{N}-|t|)^{-p/4} dt)^{1/p}
   \le 2\,c_p (\int\limits_{\RR}^{} \|f(s)\|^p
   |\sqrt{N}-|s|\; |^{-p/4} ds)^{1/p}$ .}
  \end{enumerate}
\newpage
 \item[]{Take $f := \sum\limits_{j=1}^{n+1} y_j {\chi}_j/|I_j|$
  to estimate $J_1$}
  \begin{enumerate}
  \item[$J_1$]{$\le c_4 (\int\limits_{\RR}^{}  \|\sum\limits_{j=1}^{n+1}
   y_j {\chi}_j(s) / |I_j|\;\|^p |\sqrt{N}-|s|\; |^{-p/4} ds)^{1/p} ,
   c_4 = 2\,c_p\,c_2\\
   = c_4 [\sum\limits_{j=1}^{n+1} \|y_j\|^p / |I_j|^p
   (\int\limits_{t_j}^{t_{j-1}} |\sqrt{N}-|s|\; |^{-p/4} ds)]^{1/p}\\
   \hspace*{-0.5em}\schl{c_5} [\sum\limits_{j=1}^{n+1} \|y_j\|^p / |I_j|^{p-1}
   (\sqrt{N}-|t_j|)^{-p/4}]^{1/p}\\
   = (\sum\limits_{j=1}^{n+1}\|q(t_j)e^{-t_j^2/2}\|^p/
   (n^{p/8}|{\cal H}_{n+1}(t_j)|^p |I_j|^{p-1} (\sqrt{N}-|t_j|)^{p/4}))^{1/p}\\
   \hspace*{-0.5em}\schl{}
   (\sum\limits_{j=1}^{n+1}{\mu}_j\|q(t_i)e^{-t_j^2/2}\|^p)^{1/p}$}
  \end{enumerate}
 \item[]{using lemma 1:
  $n^{1/8}|{\cal H}_{n+1}'(t_j)| |I_j| (\sqrt{N}-|t_j|)^{1/4} \schl{} 1 ,
  |I_j| \schl{} {\mu}_j$ .\\

  To estimate $J_2$, we note that there is $c_6$ such that for all
  $\ell\in{\NN}$ and\\
  $t\in (\ell\,2\sqrt{N} , (\ell +1) 2\sqrt{N})$:
  $f_j(t) \le c_4 {\ell}^{-9/4} f_j(t-2\ell\sqrt{N})$ .
  Moreover for $t\in (\sqrt{N} , 2\sqrt{N})$ and $s := 2\sqrt{N}-t
  \in (0,\sqrt{N})$ one finds  $f_j(t) \le f_j(s)$ .
  Similar statements hold for $t < -\sqrt{N}$ .
  Since $\sum\limits_{\ell}^{} {\ell}^{-9/4} < \infty$,
  this implies that there is $c_7$ such that}
  \begin{enumerate}
  \item[$J_2$]{$\le c_7 (\int\limits_{-2\sqrt{N}}^{2\sqrt{N}}
   (\sum\limits_{j=1}^{n+1}\|y_j\| f_j(t))^p dt)^{1/p} \le
   2 c_7 (\int\limits_{-\sqrt{N}}^{\sqrt{N}} (\sum\limits_{j=1}^{n+1}
   \|y_j\| f_j(t))^p dt)^{1/p}\\
   \hspace*{-0.5em}\schl{c_8}\{\sum\limits_{k=1}^{n+1} |I_k|
   (\sum\limits_{\stackrel{j=1}{j\neq k}}^{n+1}
   \|y_j\| |I_j| / (t_k-t_j)^2 + \|y_k\|/|I_k|)^p
   (\sqrt{N}-|t_k|)^{-p/4} \}^{1/p}\\
   \le (\sum\limits_{k=1}^{n+1} 
   (\sum\limits_{\stackrel{j=1}{j\neq k}}^{n+1}
   \frac{|I_j| |I_k|^{1/p}}{(t_k-t_j)^2}
   \frac{\|y_j\|}{(\sqrt{N}-|t_k|)^{1/4}})^p)^{1/p}
   + (\sum\limits_{k=1}^{n+1} \|y_k\|^p / |I_k|^{p-1}
   \cdot (\sqrt{N}-|t_k|)^{-p/4})^{1/p}\\
   =: J_{21} + J_{22}$ .}
  \end{enumerate}
 \item[]{The step involving $c_8$ is by discretization, decomposing
  $(-\sqrt{N},\sqrt{N})$ into the intervals $I_j$ and using
  the definition of $f_j$. The term $J_{22}$ is estimated as before;
  using the definition of $y_j$ and lemma 1, we find}
  \begin{enumerate}
  \item[]{$J_{21} \sim (\sum\limits_{k=1}^{n+1} | \sum\limits_{j=1}^{n+1}
   a_{kj} (\mu_j^{1/p} \|g(t_j) e^{-t_j^2/2}\|)|^p)^{1/p}$}
  \end{enumerate}
 \item[]{where  $a_{kj} = \mu_j^{3/2-1/p} \mu_k^{1/2+1/p} / (t_k-t_j)^2$
  for $j\neq k$ and $a_{kk}=0$. To bound $J_{21}$ by
  $M (\sum\limits_{j=1}^{n+1} \mu_j\|g(t_j) e^{-t_j^2/2}\|^p)^{1/p}$,
  we have to show that $A = (a_{kj})_{k,j=1}^{n+1}$ has norm $\le M$
  as a map from $\ell_p^{n+1}$ to $\ell_p^{n+1}$, $M$ being independent
  of $n\in{\NN}$.
  To do so, use lemma 3 with $\Omega = \{1,\ldots ,n+1\} , \mu \{j\} = 1 ,
  r_{kj} = (\mu_j / \mu_k )^{1/pp'}$ .
  Calculation shows that the two conditions in (12) reduce
  to one condition , namely}
  \begin{enumerate}
  \item[]{$\sup\limits_{j=1,\ldots ,n+1}
   \sum\limits_{\stackrel{k=1}{k\neq j}}^{n+1}
   \mu_k^{3/2} \mu_j^{1/2} / (t_k-t_j)^2 \le M$ .}
  \end{enumerate}
 \item[]{We check this using lemma 1:
  $\mu_k \sim |t_k-t_{k+1}| \le |t_k-t_j|$ for $k \neq j$.
  Thus we may replace discrete sums by integrals to find
  with constants independent of $j$ and $n$}
  \begin{enumerate}
  \item[]{$\sum\limits_{\stackrel{k=1}{k\neq j}}^{n+1}
   \mu_k^{3/2} \mu_j^{1/2} / (t_k-t_j)^2\\
   \le c_9 \sqrt{\mu_j} \sum\limits_{k\neq j}^{} \mu_k / |t_k-t_j|^{3/2}\\
   \le c_{10} \sqrt{\mu_j} \int\limits_{|t-t_j|\ge \mu_j}^{}
   |t-t_j|^{-3/2} dt = 4 c_{10} =: M$ .}
  \end{enumerate}
 \item[]{The estimates for $J_1$ and $J_{21},J_{22}$
  together yield inequality (3).}
 \end{enumerate}
\end{enumerate}

\begin{enumerate}
\item[(b)]{(2)
 This statement is a dualization of the inequality just proved:
 Let $4/3 < p < \infty $ and $X$ be a UMD-space.
 For any $q\in\Pi_n(X)$ there are functionals $\xi_j\in X^*$ with
 $(\sum\limits_{j=1}^{n+1} \mu_j
 \|\xi_j e^{-t_j^2/2}\|_{X^*}^{p'})^{1/p'} = 1$ such that}
 \begin{enumerate}
 \item[]{$(\sum\limits_{j=1}^{n+1} \mu_j \|q(t_j) e^{-t_j^2/2}\|^p)^{1/p}
  = \sum\limits_{j=1}^{n+1} \mu_j \; <q(t_j),\xi_j> \, e^{-t_j^2}\\
  = \sum\limits_{j=1}^{n+1} \lambda_j \; <q(t_j),\xi_j> =: I$ \ \ .}
 \end{enumerate}
\item[]{Let $r := \sum\limits_{j=1}^{n+1} \xi_j \ell_j\in\Pi_n(X^*)$ .
 Then $<q,r> \; \in \Pi_{2n}$ is integrated exactly by
 Gaussian quadrature, and H\"older's inequality yields}
 \begin{enumerate}
 \item[]{$I
  = \int\limits_{\RR}^{} <q(t),r(t)> e^{-t^2} dt\\
  \le (\int\limits_{\RR}^{} \|q(t) e^{-t^2/2}\|^p dt)^{1/p}
  (\int\limits_{\RR}^{} \|r(t) e^{-t^2/2}\|_{X^*}^{p'} dt)^{1/p'}$ .}
 \end{enumerate}
\item[]{Since $X^*$ is a UMD-space as well and $1<p'<4$, we have by (3)}
 \begin{enumerate}
 \item[]{$(\int\limits_{\RR}^{}\|r(t) e^{-t^2/2}\|_{X^*}^{p'} dt)
  ^{1/p'} \le
  c_{p'} (\sum\limits_{j=1}^{n+1}
  \mu_j \|\xi_j e^{-t_j^2/2}\|_{X^*}^{p'})^{1/p'} = c_{p'}$ .}
 \end{enumerate}
\end{enumerate}

\begin{enumerate}
\item[(b)]{(1) We denote by\\
   ${\cal K}_j(t,s) = \sum\limits_{i=0}^{j}
   {\cal H}_i(t) {\cal H}_i(s) ,\;\; {\cal K}^m(t,s) = \frac{1}{m}
   \sum\limits_{j=0}^{m-1} {\cal K}_j(t,s)$
   the kernels of the orthogonal projection $P_j$ onto
   $\Pi_j\cdot e^{-t^2/2} \le L_2({\RR})$
   and the first C\'esaro mean operator $\sigma_m$ in $L_2({\RR})$,
   respectively. It was shown by Freud $[$F1$]$ and independently by
   Poiani $[$Po$]$ that}
   \begin{enumerate}
   \item[]{$\sup\limits_{m\in{\NN}} \sup\limits_{s\in{\RR}}\;
    \int\limits_{\RR}^{} |{\cal K}^m(t,s)| dt \le M$\hfill (15)}
   \end{enumerate}
  \item[]{Let $0<\delta <1$. By Nevai $[$N,p.265$]$ there is $c_{\delta}\ge 1$
   such that for all $q\in \Pi_{4n}$}
   \begin{enumerate}
   \item[]{$\sum\limits_{|t_j|\le \delta \sqrt{N}}^{}
    \mu_j |q(t_j) e^{-t_j^2/2}|
    \le c_{\delta} \int\limits_{\RR}^{} |q(t) e^{-t^2/2}| dt$.}
   \end{enumerate}
  \item[]{We apply this to \ $q(t) e^{-t^2/2} = {\cal K}^m(t,s)$ \ for
   fixed $s\in{\RR}$ and $m\le 4n$ to get}
   \begin{enumerate}
   \item[]{$\sup\limits_{m\le 4n}  \sup\limits_{s\in{\RR}}
    \sum\limits_{|t_j|\le\delta\sqrt{N}}
    \mu_j |{\cal K}^m(t_j,s)| \le c_{\delta}\cdot M$\hfill (16)}
   \end{enumerate}
  \item[]{From (15) and (16), we find for $m\le 4n$ and $f\in L_p({\RR},X),
   \;\; p = 1$ or $\infty$ ,}
   \begin{enumerate}
   \item[]{$\sum\limits_{|t_j|\le\delta\sqrt{N}} \mu_j
    \|\sigma_m f(t_j)\| \le \sup\limits_{s}
    \sum\limits_{|t_j|\le\delta\sqrt{N}}^{} \mu_j |{\cal K}^m(t_j,s)|
    \cdot \|f\|_{L_1(X)} \le c_{\delta} M \|f\|_{L_1(X)}$\\[0.5ex]
    $\sup\limits_{|t_j|\le\delta\sqrt{N}} \|\sigma_m f(t_j)\|  \le  
    \sup\limits_{j} \int\limits_{\RR}^{} |{\cal K}^m(t_j,s)| ds\cdot
    \|f\|_{L_{\infty}(X)} \le M \|f\|_{L_{\infty}(X)}$.}
   \end{enumerate}
  \item[]{Hence, by the Riesz-Thorin interpolation theorem, for
   $f\in L_p(X) , m\le 4n$}
   \begin{enumerate}
   \item[]{$(\sum\limits_{|t_j|\le\delta\sqrt{N}}
    \mu_j \|\sigma_m f(t_j)\|^p)^{1/p}
    \le c_{\delta} M\cdot (\int\limits_{\RR} \|f(t)\|^p dt)^{1/p}$.
    \hfill (17)}
   \end{enumerate}
  \item[]{Replacing $c_{\delta} M$ by $3c_{\delta} M$, we may substitute the
   operator $v_{2n} := 2\sigma_{4n}-\sigma_{2n} =
   \frac{1}{2n} \sum\limits_{j=2n}^{4n-1} P_j$ for $\sigma_m$ in (17).
   However, $V_{2n} f = f$ for functions of the form
   $f(t) =  q(t) e^{-t^2/2}$ where $q\in \Pi_{2n}(X)$.
   Thus for any $q\in \Pi_{2n}(X)$}
   \begin{enumerate}
   \item[]{$(\sum\limits_{|t_j|\le\delta\sqrt{N}}
    \mu_j \|q(t) e^{-t^2/2}\|^p)^{1/p}
    \le 3 c_{\delta} M (\int\limits_{\RR}^{}
    \|q(t)e^{-t^2/2}\|^p dt)^{1/p}$.}
   \end{enumerate}
  \item[]{This ends the proof of (b) of theorem 1.}
\end{enumerate}

\underline{Remark.}
\begin{enumerate}
\item[]{It is likely that (5) holds for all $p$ and $X$; however, this proof
 does not work: for $m = (1+{\epsilon})n$ with ${\epsilon}>0$}
 \begin{enumerate}
 \item[]{${\mu}_1 |{\cal K}^m(t_1,t_1)| \schl{c\raisebox{-0.75ex}{$\epsilon$}}
  n^{-1/6}\cdot n^{1/2} = n^{1/3}$}
 \end{enumerate}
\item[]{tends to $\infty$ with $n$, an thus (16) does not hold if the sum
 is extended over all $j=1,\ldots ,n+1$. On the other hand,}
 \begin{enumerate}
 \item[]{$\sup\limits_{m\le n} \sup\limits_{s\in{\RR}}\;
  \sum\limits_{j=1}^{n+1} {\mu}_j |{\cal K}^m(t_j,s)| \le  M$\hfill (18)}
 \end{enumerate}
\item[]{is correct; Freud's rather elegant proof of (15) in $[$F1$]$ may be
 modified to yield (18) in the discrete case. The reason why this only works
 for $m\le n$ is that the biorthogonality relations}
 \begin{enumerate}
 \item[]{${\delta}_{kl} = \int\limits_{\RR} h_k(t)h_l(t)e^{-t^2}dt =
  \sum\limits_{j=1}^{n+1} {\lambda}_jh_k(t_j)h_l(t_j)$}
 \end{enumerate}
\item[]{are used which in the discrete case only hold if $k+l\le 2n+1$,
 i.e. essentially\\
 $k,l\le n$ is satisfied. Instead of $V_{2n}$ as above,
 one might take ${\epsilon}>0$ and\\
 $V_{\textstyle\epsilon} := \frac{1}{\textstyle\epsilon} {\sigma}_n -
 \frac{1-{\textstyle\epsilon}}{\textstyle\epsilon}
 {\sigma}_{(1-\textstyle\epsilon )n}$
 (assuming WLOG that ${\epsilon}n\in {\NN}$) to find an inequality of
 type (5) for polynomials of degree $\le (1-\epsilon )n$ and of type (3)
 by dualization without assumptions on $X$ and $p$. Thus one finds the}
\end{enumerate}
\vspace{1ex}
\underline{Proposition.}
\begin{enumerate}
\item[]{Let $1\le p\le\infty$, $X$ be a Banach space and $0<\epsilon <1$.
 Then there is $c_{\textstyle\epsilon}>0$
 such that for all $m,n\in {\NN}$ with
 $m\le (1-\epsilon )n$ and all $q\in {\Pi}_m(X)$\\[0.5ex]
 $c_{\textstyle\epsilon}^{-1}
 (\sum\limits_{i=1}^{n+1} {\mu}_i \|q(t_i)e^{-t_i^2/2}\|^p)^{1/p} \le
 (\int\limits_{\RR} \|q(t)e^{-t^2/2}\|^pdt)^{1/p} \le
 c_{\textstyle\epsilon}
 (\sum\limits_{i=1}^{n+1} {\mu}_i \|q(t_i)e^{-t_i^2/2}\|^p)^{1/p}$.}
\end{enumerate}

This should be compared with theorem 1. We now return to the one part of
theorem 1 still to be proved.

\begin{enumerate}
\item[(a)]{$(1)\Rightarrow (2)$.}
 \begin{enumerate}
 \item[]{Assume (1) holds. We have to show that necessarily
  $1<p<4$ holds and that $X$ is a UMD-space.\\
  Choose $q=h_n$ in inequality (3). Since by (2),
  $|{\cal H}_n(t_j)|=\sqrt{2/{(n{\mu}_j)}}$, the right side in (3) reads}
  \begin{enumerate}
  \item[]{$(\sum\limits_{j=1}^{n+1} {\mu}_j |{\cal H}_n(t_j)|^p)^{1/p} = 
   \sqrt{\frac{2}{n}} (\sum\limits_{j=1}^{n+1} {\mu}_j^{1-p/2})^{1/p} \le
   \sqrt{\frac{2}{n}}(\sum\limits_{j=1}^{n+1} {\mu}_j)^{1/p}
   \sup\limits_{j} {\mu}_j^{-1/2} \\[0.5ex]
   \schl{}  n^{-1/2}n^{1/{2p}}n^{1/4} = n^{1/{2p}-1/4}$}
  \end{enumerate}
 \item[]{by using lemma 1 again. For the left side, we find using only the
  asymptotic behavior of ${\cal H}_n$ near its maximum $[$AW$]$}
  \begin{enumerate}
  \item[]{$\|{\cal H}_n\|_p \gschl n^{-1/{6p}-1/{12}}$ \ .}
  \end{enumerate}
\newpage
 \item[]{Hence (3) requires $-\frac{1}{6p} -\frac{1}{12} \le
  \frac{1}{2p} -\frac{1}{4}$, i.e. $p\le 4$. For $p=4$, a slightly more
  careful use of the formulas in $[$AW$]$ shows}
  \begin{enumerate}
  \item[]{$\|{\cal H}_n\|_4 \schl{} n^{-1/8}(\log n)^{1/4}$}
  \end{enumerate}
 \item[]{which is larger than $n^{-1/8}$ and excludes $p=4$ as well.\\
 
  We now show that (3) implies that the Hilbert matrix
  $A=(1/(i-j+1/2))_{i,j\in {\NN}}$ defines a bounded operator
  $A:l_p(X)\longrightarrow l_p(X)$.
  Since $A$ is not bounded in $l_1$, this excludes $p=1$.
  A well-known approximation and scaling argument shows that the boundedness
  of $A$ in $l_p(X)$ is equivalent to the boundedness of the Hilbert transform
  in $L_p({\RR};X)$, i.e. $X$ is a UMD-space.
  In this sense, $A$ is a discrete version of the Hilbert transform.
  For $n\in {\NN}$, we need the zeros $(t_j^{n+1})$ of $H_{n+1}$ and
  $(t_i^n)$ of $H_n$. Let
  $J:=\{j\;|\;\;|t_j^{n+1}|\le 1\}$, $I:=\{i\;|\;\;|t_i^n|\le 1\}$
  . Take any system
  $(x_j)_{j\in J}\subseteq X$ and define $q\in {\Pi}_n$ by}
  \begin{enumerate}
  \item[]{$q(t) :=
   \sum\limits_{j\in J} 
   ({\mu}_j^{n+1})^{-1/p}x_je^{(t_j^{n+1})^2/2}l_j^{n+1}(t)$,}
  \end{enumerate}
 \item[]{where $l_i^{n+1}\in{\Pi}_n, l_i^{n+1}(t_k^{n+1})={\delta}_{jk}$
  for $i,k=1,\ldots ,n+1$. Note that $q(t_k)=0$ for $k\not\in J$. Thus using
  assumption (3) and inequality (4) of (b)(1) of
  theorem \nolinebreak 1 -- but with
  ${\mu}_i^n, t_i^n$ instead of ${\mu}_j^{n+1}, t_j^{n+1}$ -- we find}
  \begin{enumerate}
  \item[]{$\sum\limits_{i\in I}
   {\mu}_i^n \|q(t_i^n)e^{-(t_i^n)^2/2}\|^p\\[0.75ex]
   \le c_1 \int\limits_{\RR} \|q(t)e^{-t^2/2}\|^p dt\\[0.25ex]
   \le c_2 \sum\limits_{j\in J}
   {\mu}_j^{n+1} \|q(t_j^{n+1})e^{-(t_j^{n+1})^2/2}\|^p$,}
  \end{enumerate}
 \item[]{i.e.}
  \begin{enumerate}
  \item[]{$\sum\limits_{i\in I} \| \sum\limits_{j\in J} a_{ij}x_j\|^p \le
   c_2 \sum\limits_{j\in J} \|x_j\|^p$.}
  \end{enumerate}
 \item[]{where}
  \begin{enumerate}
  \item[]{$a_{ij} :=
   ({\mu}_i^n/{\mu}_j^{n+1})^{1/p}l_j^{n+1}(t_i^n)
   e^{(t_j^{n+1})^2/2}e^{-(t_i^n)^2/2}\\[0.25ex]
   = ({\mu}_i^n/{\mu}_j^{n+1})^{1/p}{\cal H}_{n+1}(t_i^n) /
   [{\cal H}_{n+1}'(t_j^{n+1})(t_i^n-t_j^{n+1})]$ for $i\in I, j\in J$.}
  \end{enumerate}
 \item[]{By (2), $1/|{\cal H}_{n+1}'(t_j^{n+1})| = \sqrt{{\mu}_j^{n+1}/2}$.
  Using the recursive formulas of the $H_n$'s, see Szeg\"o $[$Sz, 5.5$]$,
  ${\cal H}_n(t_i^n)=0$, and again (2), one finds that}
\newpage
  \begin{enumerate}
  \item[]{$|{\cal H}_{n+1}(t_i^n)| = |-\sqrt{n/(n+1)}{\cal H}_{n-1}(t_i^n)|
   = \sqrt{n/(n^2-1)}/\sqrt{{\mu}_i^n}$,\\[0.25ex]
   $a_{ij} = {\epsilon}_i{\delta}_j({\mu}_i^n/{\mu}_j^{n+1})^{1/p-1/2}
   \sqrt{n/(2(n^2-1))}/(t_i^n-t_j^{n+1})$,}
  \end{enumerate}
 \item[]{${\epsilon}_i = \sgn\;{\cal H}_{n+1}(t_i^n),
  {\delta}_j = \sgn\;{\cal H}_{n+1}'(t_j^{n+1})$. By lemma 1,
  ${\mu}_i^n \schl{c_3} \frac{1}{n} \schl{c_3} {\mu}_j^{n+1}$ for $i\in I$,
  $j\in J$. Hence $B = (b_{ij})$ with
  $b_{ij} = \frac{1}{\sqrt{N} (t_i^n-t_j^{n+1})},\;\; N = 2n+3$ defines a map\\
  $B: l_p^{|J|}(X)\longrightarrow l_p^{|I|}(X)$ of norm $\le b_2b_3$
  bounded independently of $n\in {\NN}$. Near zero, the asymptotic formula}
  \begin{enumerate}
  \item[]{${\cal H}_{n+1}(t) = \frac{1}{\sqrt{\pi}}(\frac{2}{n})^{1/4}
   [\cos{(\sqrt{N}\,t-n\frac{\pi}{2} )} +
   \frac{t^3}{6\sqrt{N}}\sin{(\sqrt{N} t-n\frac{\pi}{2}) +
   {\cal O} (\frac{1}{n})}],\\[0.25ex]
   N = 2n+3$}
  \end{enumerate}
 \item[]{for the Hermite functions holds, $[$Sz, 8.22.6$]$. The zeros of
  ${\cal H}_{n+1}$ in $[-1,1]$ may be determined from the $\cos$-term up to
  an error of ${\cal O} (\frac{1}{n})$, since
  $(\sqrt{N} t_j^{n+1}-n\frac{\pi}{2} )$ is determined up to
  ${\cal O} (\frac{1}{\sqrt{N}} )$. The zeros of ${\cal H}_{n+1}$ separate
  those of ${\cal H}_n$; for ${\cal H}_n$, $N=2n+3$ is replaced by
  $\tilde{N} = 2n+1$ with (again) $\sqrt{N} -\sqrt{\tilde{N}} =
  {\cal O} (\frac{1}{\sqrt{N}} )$. Thus the difference $t_j^n-t_j^{n+1}$ is
  $\frac{\pi}{2\sqrt{N}}+{\cal O} (\frac{1}{n} )$ and}
  \begin{enumerate}
  \item[]{$|(t_i^n-t_j^{n+1}) - \frac{\pi (i-j+1/2)}{\sqrt{N}} | \le
   c_4/n$; $i\in I$, $j\in J$}
  \end{enumerate}
 \item[]{and thus}
  \begin{enumerate}
  \item[]{$|b_{ij}-\frac{1}{\pi (i-j+1/2)} | \le c_5 \frac{1}{(i-j+1/2)^2} =:
   c_{ij}$ ; \ $i\in I$, $j\in J$}
  \end{enumerate}  
 \item[]{where $c_4$, $c_5$ are independent of $n$, $i\in I$ and $j\in J$.
  Since $C = (c_{ij} )_{i,j\in {\NN}}$ is bounded as a map
  $l_p(X)\longrightarrow l_p(X)$, we conclude that the finite Hilbert matrix
  $(\frac{1}{i-j+1/2} )_{i\in I, j\in J} :
  l_p^{|J|}(X)\longrightarrow l_p^{|I|}(X)$ has norm independent of
  $n\in {\NN}$, i.e. of $I$, $J$. Note that $|I| \schl{}  |J| \schl{}
  \sqrt{n} \longrightarrow \infty$ for $n\longrightarrow\infty$. Thus $A$
  is continuous in $l_p(X)$.
  This ends the proof of theorem 1.\hfill\raisebox{-2ex}{$\Box$}\\}
 \end{enumerate}
\end{enumerate}

\underline{Remarks.}
\begin{enumerate}
\item[(a)]{The proof shows that $X$ has to be a UMD-space already if (3)
 only holds for those $q\in\Pi _n(X)$ with $q(t_j) = 0$ for all $|t_j|>1$.}
\item[(b)]{Theorem 1 shows that
 $(\int\limits_{\RR} |q(t)e^{-t^2/2}|^pdt)^{1/p}$ is essentially determined
 by the values of $q$ in $(-\sqrt{N} ,\sqrt{N} )$, see also Freud $[$F2$]$.
 This corresponds to the fact that, if $|q(\bar{t})|e^{-{\bar{t}}^2/2}
 = \max\limits_{t\in{\RR}} |q(t)| e^{-t^2/2}$
 for $q\in\Pi _n$, then $|\bar{t} |\le\sqrt{N}$, as can be shown using Gaussian
 quadrature techniques.}
\item[(c)]{For a UMD-space and $4/3<p<4$, the subspace $\Pi _n(X)$ of
 $L_p({\RR} ,e^{-t^2/2p};X)$ is uniformly isomorphic to $l_p^{n+1}(X)$,
 the maps}
 \begin{enumerate}
 \item[]{$J_n: \Pi _n(X)\longrightarrow l_p^{n+1}(X),\\
  q\longrightarrow ({\mu}_j^{1/p}e^{-t_j^2/2}q(t_j))_{j=1}^{n+1}$}
 \end{enumerate} 
\item[]{satisfy $\sup\limits_{n} \|J_n\|\;\|J_n^{-1}\|\le c_p<\infty$,
 interpolation essentially yields the Banach-Mazur distance. This probably
 also holds for $1<p\le 4/3$.\\}
\end{enumerate}

\section{Mean convergence of interpolating polynomials and expansions}

We now give the \underline{Proof of theorem 2} :

Let X be a UMD-space, $1<p<4$ and $\alpha >1/p$. Define the norm
\begin{enumerate}
\item[]{$\dline g\dline :=\sup_{t\in {\RR}}
 \|g(t)\| _X(1+|t|)^{\alpha}e^{-t^2/2}$}
\end{enumerate}
for those $g\in C({\RR} ;X)$ where this is finite. Take
$f\in C({\RR} ;X)$ satisfying (6). Then $\dline f\dline <\infty$, and
moreover,
$f$ can be approximated by polynomials $q_n\in \Pi _n(X)$ in
$\dline\cdot\dline$,
$\dline f-q_n\dline\longrightarrow 0$ ((6) allows the restriction to a finite
intervall where this clearly is possible).\\
Let $\| g\| _p:= (\int\limits_{\RR} \| g(t)\| ^pe^{-t^2/2p}dt)^{1/p}$.
Since $q_n-I_nf\in\Pi _n(X)$ and $I_nf(t_j)=f(t_j)$, theorem 1 yields
\begin{enumerate}
\item[]{$\| f-I_nf\| _p\le\| f-q_n\| _p+\| q_n-I_nf\| _p\\
 \le (\int\limits_{\RR}
 \frac{dt}{(1+|t|)^{\alpha p}})^{1/p}\dline f-q_n\dline
 +(\sum\limits_{j=1}^{n+1} {\mu}_j\| (q_n(t_j)-f(t_j))e^{-t_j^2/2}\| ^p)^{1/p}$}
\end{enumerate}
The integral is finite $(\alpha p>1)$, thus the first term tends to zero.
The second term approaches zero just as well: using again
$\mu _j\schl{} (t_j-t_{j+1})$, see lemma 1, we find
\begin{enumerate}
\item[]{$(\sum\limits_{j=1}^{n+1}
 \mu _j\| (f(t_j)-q_n(t_j))e^{-t_j^2/2}\| ^p)^{1/p} \le
 \dline f-q_n\dline\cdot (\sum\limits_{j=1}^{n+1}
 \mu _j/(1+|t_j|)^{\alpha p})^{1/p}\\[0.25ex]
 \le
 c\dline f-q_n\dline\cdot (\int\limits_{\RR}dt/(1+|t|)^{\alpha p})^{1/p}$.}
\end{enumerate}
Hence $\| f-I_nf\| _p\longrightarrow 0$.

For $p>4$ the same statement does not hold, as we will show now.
The example is an extension and modification of Nevai's $[$N$]$. Assume $p>4$
and choose $\alpha$ with $1/p<\alpha <1/4$. Consider the Banach spaces
\begin{enumerate}
\item[]{$C_0:= (\{f: {\RR}\longrightarrow {\RR}\;|\;\,
 f$ continuous, $|f(t)|(1+|t|)^{\alpha} e^{-t^2/2}\longrightarrow 0$ for $
 |t|\longrightarrow\infty\}, \dline\cdot\dline )$\\
 $L_p:= \{f: {\RR}\longrightarrow {\RR}\;|\;\;
 \| f\| _p= (\int\limits_{\RR} (|f(t)|e^{-t^2/2})dt)^{1/p}<\infty \}$.}
\end{enumerate}
If theorem 2 would hold for this $\alpha$ and $p$, the interpolating operators
$I_n: C_0\longrightarrow L_p,\; f\longmapsto I_nf$ would be uniformly
bounded by the Banach-Steinhaus theorem,\\
$\sup\limits_{n\in {\NN}}\| I_n: C_0\longrightarrow L_p\|=M<\infty$.\\[2ex]
Let $\epsilon _j:= \sgn\;{\cal H}_{n+1}'(t_j),\;J:= \{j\,|\, t_j\le 0\}$
and define
\begin{enumerate}
\item[]{\[ f(t_j):= \left\{ \begin{array}{c@{\quad:\quad}l}
 \epsilon _j(1+|t_j|)^{-\alpha }e^{t_j^2/2} & j\in J \\
 0 & j\not\in J \end{array} \right\}\;\;\;\raisebox{-2ex}{.} \] }
\end{enumerate}

\vspace{2ex}
This obviously can be extended to define a continuous function
$f:{\RR}\longrightarrow {\RR}$ with $\dline f\dline =1$.
We show that $I_nf$
is bounded from below by a suitable multiple of $h_{n+1}$, for
$0\le t\le\sqrt{N}$: Then $0\le t-t_j\le 2\sqrt{N}$ for $j\in J$. Using (2),
$\mu _j\ge c_1/\sqrt{n}$ (lemma 1) and $|J|\schl{} n/2$, we find
\begin{enumerate}
\item[]{$|I_nf(t)|e^{-t^2/2}= \sum\limits_{j\in J}
 \frac{|{\cal H}_{n+1}(t)|}{|{\cal H}_{n+1}'(t_j)|(1+|t_j|)^{\alpha }(t-t_j)}
 \ge c_2 \sum\limits_{j\in J}
 \sqrt{{\mu}_j} n^{-(1+\alpha )/2}|{\cal H}_{n+1}(t)|\\[0.25ex]
 \ge c_3 n^{-1/4-\alpha /2}|{\cal H}_{n+1}(t)|,
 \hspace{5cm} 0\le t\le\sqrt{N}$.}
\end{enumerate}
Using Skovgaard's formula (8), $p>4$ and $\alpha <1/4$, we get for all
$n\in {\NN}$
\begin{enumerate}
\item[]{$M\ge \| I_n\| \dline f\dline \ge \| I_nf\| _p \ge
 c_3\,n^{1/4-\alpha /2}(\int\limits_{0}^{\sqrt{N}}
 |{\cal H}_{n+1}(t)|^pdt)^{1/p} \schl{}
 n^{1/4-\alpha /2}n^{-1/12-1/{6p}}\\[0.25ex]
 = n^{1/6-1/{6p}-\alpha /2} \ge
 n^{1/24-1/{6p}}$ \ .}
\end{enumerate}
Thus $1/24-1/{6p}\le 0,\; p\le 4$, contradicting our assumption $p>4$.\\
This proves theorem 2.\hfill\raisebox{-2ex}{$\Box$}\\[2ex]

In a similar way we may prove that the assumption
\begin{enumerate}
\item[]{$f\in C({\RR}) \cap L_p({\RR} ;e^{-t^2/2}) =:Y$}
\end{enumerate}
does not suffice, in general, to prove the convergence $I_nf\longrightarrow f$
in $L_p$-norm (Remark (ii) after theorem 2):\\
Assume this would hold. Introduce the norm
\begin{enumerate}
\item[]{$\| f\| _0:= \max (\sup\limits_{t\in {\RR}} |f(t)|e^{-t^2/2},
 (\int\limits_{\RR} (|f(t)|e^{-t^2/2})^pdt)^{1/p}), f\in Y$}
\end{enumerate}
on $Y$. Then $(Y,\|\cdot\| _0)$ is a Banach space and by Banach-Steinhaus
we would have\\
$\| I_nf\| _p\le M\| f\| _0$. Define $f\in Y$ as above, with
$\alpha =0$, and $\| f\| _0=1$. By the above estimates
(for all $1\le p<\infty , n\in {\NN}$)
\begin{enumerate}
\item[]{$M\ge \| I_nf\| _p\ge c_4 n^{1/4}
 (\int\limits_{0}^{\sqrt{N}} |{\cal H} _{n+1}(t)|^pdt)^{1/p}$.}
\end{enumerate}
The right side is of order $n^{1/{6p'}}$ for $p>4$, $n^{1/8}(\log n)^{1/4}$
for $p=4$ and $n^{1/{2p}}$ for $p<4$. In any of these cases, this contradicts
the uniform boundedness by $M$.
\hfill\raisebox{-3ex}{$\Box$}\\[1ex]

It remains to give the \underline{Proof of theorem 3} :

Recall that we put $a_j= \int\limits_{\RR} f(t){\cal H} _j(t)dt,\;
P_nf= \sum\limits_{j=0}^{n} a_j{\cal H} _j$ for $f\in L_p({\RR} ;X)$.
\vspace{1ex}
\begin{enumerate}
\item[$(2)\Rightarrow (1)$]{If $X$ is a UMD-space and $4/3<p<4,
 P_nf\longrightarrow f$ in $L_p({\RR} ;X)$. The scalar proof of Askey-Wainger
 $X={\RR}$ $[$AW$]$ directy generalizes to the $X$-valued case, by using the
 boundedness of the Hilbert-transform and of the weighted Hilbert-transforms
 with kernels $\frac{1}{t-s} |\frac{t}{s} |^{\pm 1/4}$ also in the vector valued
 case of $L_p({\RR} ;X)$ for $4/3<p<4$ (the latter folows using lemma 4).}
\item[$(1)\Rightarrow (2)$]{We assume that $P_nf\longrightarrow f$ for all
 $f\in  L_p({\RR} ;X)$. Already in the scalar case, $4/3<p<4$ is necessary
 $[$AW$]$. By the Banach-Steinhaus theorem,\\
 $\sup\limits_{n\in {\NN}}
 \|P_n:L_p({\RR} ;X)\longrightarrow L_p({\RR} ;X)\|= M<\infty$.
\newpage
 Using this, we can dualize inequality (4). Let $0<\delta <1$.
 We claim that there is $a_{\delta}$ such that for all $q\in\Pi _n(X)$ with
 $q(t_j)=0$ for $|t_j|>\delta\sqrt{N}$, one has}
 \begin{enumerate}
 \item[]{$(\int\limits_{\RR} \| q(t)e^{-t^2/2}\| ^pdt)^{1/p} \le
  a_{\delta} (\sum\limits_{j=1}^{n+1}
  \mu _j\|q(t_j)e^{-t_j^2/2}\| ^p)^{1/p}$\hfill (19)}
 \end{enumerate}
\item[]{By remark (a) at the end of section 3, this will imply that $X$
 has to be a UMD-space. For such $q$, let $f(t):= q(t)e^{-t^2/2}$.
 Since $P_n$ projects onto $\Pi _n(X)\cdot e^{-t^2/2}$, $P_nf=f$. Let
 $\epsilon >0$.
 Then there is $g\in L_p({\RR} ;X^{*}), \int\limits_{\RR}
 \| g(t)\| _{X^{*}}^{p'} =1$, such that}
 \begin{enumerate}
 \item[]{$(1-\epsilon )(\int\limits_{\RR}
  \| q(t)e^{-t^2/2}\| _X^pdt)^{1/p} = (1-\epsilon )(\int\limits_{\RR}
  \| f(t)\| _X^pdt)^{1/p} \le
  \int\limits_{\RR} <g(t),f(t)> dt\\[0.5ex]
  = \int\limits_{\RR} <g(t),P_nf(t)> dt
  = \int\limits_{\RR} <P_ng(t),f(t)> dt =: I$}
 \end{enumerate}
\item[]{Since $P_ng(t)$ is of the form $r(t)e^{-t^2/2}$ for some
 $r\in\Pi _n(X^{*}),\\
 <P_ng(t),f(t)> = <r(t),q(t)>e^{-t^2}$ is integrated
 exactly by Gaussian quadrature. Using $q(t_j)=0$ for $|t_j|>\delta\sqrt{N}$,
 inequality (4) for $X^{*}$ (being UMD as well) and $4/3<p'<4$, we find}

 \begin{enumerate}
 \item[$I$]{$= \sum\limits_{|t_j|\le\delta\sqrt{N}}
  \lambda _j<r(t_j),q(t_j)>\\[0.5ex]
  = \sum\limits_{|t_j|\le\delta\sqrt{N}} \mu _j<P_ng(t_j),f(t_j)>\\[0.5ex]
  \le (\sum\limits_{|t_j|\le\delta\sqrt{N}}
  \mu _j \| P_ng(t_j)\| _{X^{*}}^{p'})^{1/{p'}}(\sum\limits_{j}
  \mu _j \| f(t_j)\| ^p)^{1/p}\\[0.5ex]
  \le c_{\delta} (\int\limits_{\RR} \| P_ng(t)\| ^{p'}dt)^{1/{p'}}
  (\sum\limits_{j} \mu _j\| q(t_j)e^{-t_j^2/2}\| ^p)^{1/p}\\[0.5ex]
  \le c_{\delta} M(\sum\limits_{j} \mu _j \|q(t_j)e^{-t_j^2/2}\|^p)^{1/p}$}
 \end{enumerate}
\item[]{which gives (19). In the last step $\|P_n=P_n^{*}:
 L_{p'}({\RR} ;X^*)\longrightarrow L_{p'}({\RR} ;X^*)\|\le M$ was used.
 This ends the proof of theorem 3.}
\end{enumerate}
\hfill$\Box$\\

\newpage
\underline{References.}

\setcounter{romzahl}{2}

\begin{tabular}{ll}

$[$A$]$ & R. Askey; Mean convergence of orthogonal series and
  Lagrange interpolation,\\
  & Acta Math.~Sci.~Hungar.~23(1972), 71-85.\\

$[$AW$]$ & R. Askey, S. Wainger; Mean convergence of expansions in
  Laguerre and\\
  & Hermite series, Amer.~J. Math.~87(1965), 695-708.\\
           
$[$BMP$]$ & A.I. Benedek, E.R. Murphy, R. Panzone; Cuestions del analisis de
          Fourier,\\
        & Notas de Algebra y Analisis 5(1974), Univ.~Bahia Blanca, Argentinia.\\

$[$FG$]$ & D.L. Fernandez, J.B. Garcia; Interpolation of Orlicz-valued
           function spaces\\
         & and U.M.D. property, Studia Math.~99(1991), 23-39.\\
           
$[$F1$]$ & G. Freud; šber die $(C,1)$- Summen der Entwicklungen nach
           orthogonalen\\
         & Polynomen, Acta Math.~Sci.~Hungar.~14(1963), 197-208.\\
           
$[$F2$]$ & G. Freud; On an inequality of Markov type,
           Soviet Math.~Dokl.~12(1971), 570-573.\\

$[$GR$]$ & J. Garcia-Cuerva, J.L. Rubio de Francia; Weighted norm
           inequalities and\\
         & related topies, North Holland, 1985.\\

$[$KN$]$ & H. K\"onig, N.J. Nielsen; Vector-valued $L_p$-convergence of
           orthogonal series\\
         & and Lagrange interpolation, preprint, 1992.\\
           
$[$N$]$  & P.G. Navai; Mean convergence of Lagrange interpolation
           \Roman{romzahl},\\
         & J. Appr.~Th.~30(1980), 263-276.\\
          
$[$P$]$  & H. Pollard; Mean convergence of orthogonal series
           \Roman{romzahl},\\
         & Transact.~AMS 63(1948), 355-367.\\
          
$[$PO$]$ & E.L. Poiani; Mean Cesaro summability of Laguerre and
           Hermite series,\\
         & Transactions AMS 173(1972), 1-31.\\
           
$[$S$]$  & J. Schwartz; A remark on inequalities of Calderon-Zygmund
           type for\\
         & vector-valued functions,
           Comm.~Pures Appl.~Math.~14(1961), 785-799.\\
          
$[$Sz$]$ & G. Szeg\"o; Orthogonal polynomials, AMS, Providence, 1959.\\

$[$Z$]$  & A. Zygmund; Trigonometric series, Cambridge Univ.~Press, 1968.\\

\end{tabular}

\end{document}